\documentclass[a4paper,10pt,leqno,twoside]{article}

\usepackage[english]{babel} \usepackage{inputenc, amsmath, amssymb , latexsym,
  epic, epsfig, rotating, fancyheadings, amsthm, pifont, empheq}

\newtheorem{definition}{Definition}[section]
\newtheorem{thm}[definition]{Theorem}

\newtheorem{lem}[definition]{Lemma}
\newtheorem{cor}[definition]{Corollary}
\newtheorem{prop}[definition]{Proposition}

\newtheorem{bem}[definition]{Remark}

\newtheorem*{propa}{Proposition A}

\newtheorem*{thmc}{Theorem C}
\newtheorem*{thmc'}{Theorem C'}
\newtheorem*{thmd}{Theorem D}

\newtheorem*{corb}{Corollary B}

\newtheorem*{core}{Corollary E}

\numberwithin{equation}{section}

\newcommand{\homtwo}{\ensuremath{\mathrm{Homeo}_0(\mathbb{T}^2)}}

\newcommand{\homd}{\ensuremath{\mathrm{Homeo}_0(\mathbb{T}^d)}}
\newcommand{\equi}{\ensuremath{\Leftrightarrow}}

\newcommand{\ld}{\ensuremath{,\ldots,}}
\newcommand{\ssq}{\ensuremath{\subseteq}}
\newcommand{\smin}{\ensuremath{\setminus}}

\newcommand{\T}{\ensuremath{\mathbb{T}}}

\newcommand{\N}{\ensuremath{\mathbb{N}}} 
\newcommand{\R}{\ensuremath{\mathbb{R}}}
\newcommand{\Z}{\ensuremath{\mathbb{Z}}}
\newcommand{\Q}{\ensuremath{\mathbb{Q}}}
\newcommand{\C}{\ensuremath{\mathbb{C}}}

\newcommand{\diam}{\ensuremath{\mathrm{diam}}}
\newcommand{\inte}{\ensuremath{\mathrm{int}}}

\newcommand{\conv}{\ensuremath{\mathrm{Conv}}}
\newcommand{\kreis}{\ensuremath{\mathbb{T}^{1}}}

\newcommand{\torus}{\ensuremath{\mathbb{T}^2}}
 
\newcommand{\alphlist}{\begin{list}{(\alph{enumi})}{\usecounter{enumi}}}
\newcommand{\romanlist}{\begin{list}{(\roman{enumi})}{\usecounter{enumi}}}
\newcommand{\listend}{\end{list}}

\newcommand{\solidqed}{{\raggedleft \large $\blacksquare$ \\}}

\newcommand{\foot}{\footnote}
\newcommand{\nLim}{\ensuremath{\lim_{n\rightarrow\infty}}}
\newcommand{\iLim}{\ensuremath{\lim_{i\rightarrow\infty}}}

\title{\Large\textsc{Linearization of conservative toral homeomorphisms}}
\author{T.~J\"ager\thanks{Coll\`ege de France, Paris. Email: {\tt
      tobias.jager@college-de-france.fr}}}

\newcommand{\A}{\ensuremath{\mathbb{A}}}
\newcommand{\homa}{\ensuremath{\mathrm{Homeo}_0(\mathbb{A})}}
\newcommand{\homeo}{\ensuremath{\mathrm{Homeo}}}
\newcommand{\homar}{\ensuremath{\mathrm{Homeo}_0^\mathrm{nw}(\mathbb{A})}}

\begin{document}

\setlength{\oddsidemargin}{0.02\textwidth}
\setlength{\evensidemargin}{0.02\textwidth}

\maketitle 

\abstract{We give an equivalent condition for the existence of a semi-conjugacy
  to an irrational rotation for conservative homeomorphisms of the
  two-torus. This leads to an analogue of Poincar\'e's classification of circle
  homeomorphisms for conservative toral homeomorphisms with unique rotation
  vector and a certain bounded mean motion property. For minimal toral
  homeomorphisms, the result extends to arbitrary dimensions. Further, we
  provide a basic classification for the dynamics of toral homeomorphisms with
  all points non-wandering.
}

\section{Introduction}

One of the earliest, and still one of the most elegant, results
in dynamical systems is Henri Poincar\'e's celebrated classification of the
dynamics of circle homeomorphisms \cite{poincare:1885}.
\begin{quote}
  \em An orientation-preserving homeomorphism of the circle is semi-conjugate to
  an irrational rotation if and only if its rotation number is irrational, and
  if only if it has no periodic orbits.
\end{quote}
Ever since, the question of linearization has been one of the central themes of
the subject -- when can the dynamics of a given system be related to those of a
linear model, as for example periodic or quasiperiodic motion on a torus?  It
seems natural to attempt to generalise Poincar\'e's result to higher
dimensions. However, so far no results in this direction exist. Partly, this is
explained by the fact that even on the two-torus, the situation which is best
understood, obstructions to linearization other than the existence of periodic
orbits appear. First of all, there does not have to be a uniquely defined
rotation vector. Instead, it is only possible in general to define a rotation
set, which is a compact convex subset of the plane
\cite{misiurewicz/ziemian:1989} (see also~(\ref{eq:3}) below for the
definition). Further, even when this rotation set is reduced to a single,
totally irrational rotation vector, a toral homeo- or diffeomorphism may have
dynamics which are very different from quasiperiodic ones, for example it can
exhibit weak mixing \cite{fayad:2002}. This is even true for toral flows. One
way to bypass these problems is to use higher smoothness assumptions on the
system, together with arithmetic conditions on the rotation vector, in order to
guarantee the existence of a smooth conjugacy. This is the content of
KAM-theory. However, in dimension greater than one, the price one has to pay for
this is to restrict to perturbative results, meaning that the considered
toral diffeomorphisms have to be close to the irrational rotation.

Here, we pursue a different direction. We show that whether or not a
conservative%
\foot{By conservative, we mean that there exists an invariant probability
  measure of full topological support. Due to the Oxtoby-Ulam Theorem, we can
  always assume that this measure is the Lebesgue measure on $\torus$, but we
  will actually not make use of this fact.}
toral homeomorphism is (topologically) semi-conjugate to an irrational rotation
is completely determined by the convergence properties of the rotation
vector. The method is inspired by the one in \cite{jaeger/stark:2006}, where an
analogous result is given for skew products over irrational rotations. However,
in order to overcome the lack of a fibred structure, a quite different
implementation of the ideas is required. The fact that it is possible to carry
these concepts over to the non-fibred setting also provides a new approach to
study the dynamics of periodic point free toral homeomorphisms (see Theorem~D
below and also \cite{jaeger:2008b}), which are not yet very well understood in
general (see \cite{kwapisz:2002,kwapisz:2003,beguin/crovisier/leroux:2004} for
some previous results and \cite{franks/misiurewicz:1990} for the statement of
the related Franks-Misiurewicz Conjecture).
\smallskip

Denote by \homd\ the class of homeomorphisms of the $d$-dimensional torus which
are homotopic to the identity. We say $f\in\homd$ is an {\em irrational
  pseudo-rotation}, if there exists a totally irrational vector $\rho \in \R^d$
and a lift $F:\R^d\to\R^d$ of $f$, such that for all $z\in\R^d$ there holds
\begin{equation}
  \label{eq:1}
  \nLim (F^n(z) - z)/n \ = \ \rho \ . 
\end{equation}
Similarly, when $K\ssq\T$ is an invariant subset and (\ref{eq:1}) holds for all
$z\in K$, then we say $f$ is an irrational pseudo-rotation on $K$. 

If $f$ is semi-conjugate to the irrational rotation $R_\rho : z \mapsto z+\rho
\bmod 1$, then it is further evident that there must be a certain rate of
convergence in (\ref{eq:1}), namely an a priori error estimate of $c/n$, for
some constant $c$ independent of $z$. In order to reformulate this, let
\begin{equation}
  \label{eq:2}
  D(n,z) \ := \ F^n(z)-z-n\rho \ . 
\end{equation}
We say an irrational pseudo-rotation $f$ (on an invariant set $K\ssq \T^d$) has
{\em bounded mean motion}, with constant $c\geq 0$ (on $K$), if there holds
$\|D(n,z)\| \leq c$ for all $n\in\Z$ and $z\in\R^d$ ($z\in K$). 

Now, it is a natural question to ask whether these two obvious necessary
conditions are already sufficient in order to guarantee the existence of a
semi-conjugacy. This is not true in general, counter-examples are given in
\cite{jaeger:2008b}. However, these examples exhibit wandering open sets, such
that one can still hope to obtain a positive result under additional recurrence
assumptions on the system. A first, quite elementary observation is the
following.
\begin{propa}
  \label{t.d-dim-linearization} Let $f\in\homd$, and suppose $K\ssq \T^d$ is a
  minimal set and $f$ is an irrational pseudo-rotation with bounded mean motion
  on $K$. Then $f_{|K}$ is regularly semi-conjugate to the irrational rotation
  on $\T^d$.%
  \foot{See Section~\ref{Minimal} for the definition of a {\em regular
      semi-conjugacy}. When $K=\T^d$, this just means that the semi-conjugacy is
    homotopic to the identity (and therefore, in particular, preserves the
    rotation vector and the bounded mean motion property).}
\end{propa}
In particular, when $f$ has bounded mean motion on all of $\T^d$, then its
restriction to any minimal subset is semi-conjugate to $R_\rho$. The analogue
statement holds for toral flows. 


The possibility of restricting to minimal subsets in Proposition~A is
particularly interesting in dimension two, since it can be combined with an old
result by Misiurewicz and Ziemian \cite{misiurewicz/ziemian:1991} in order to
obtain the following consequence.
\begin{corb}
  Suppose the rotation set of $f\in\homtwo$ has non-empty interior. Then for any
  totally irrational vector $\rho$ in the interior of the rotation set, there
  exists a minimal subset $K_\rho$, such that $f_{|K_\rho}$ is regularly
  semi-conjugate to $R_\rho$.
\end{corb}
This can be seen as a natural analogue of the fact that rational rotation
vectors in the interior of the rotation set are realised by periodic orbits
\cite{franks:1989}. \smallskip

In order to obtain an analogous result for conservative homeomorphisms of the
two-torus, an important ingredient will be the concept of a {\em circloid},
which is a subset $C\ssq\torus$ which is (i) compact and connected, (ii) {\em
  essential} (not contained in any embedded topological disk), (iii) has a
connected complement which contains an essential simple closed curve and (iv)
does not contain any strictly smaller subset with properties (i)--(iii).  The
semi-conjugacy in the conservative case will be obtained by constructing a
``lamination'' on the torus consisting of pairwise disjoint circloids, on which
$f$ acts in the same way as the irrational rotation on the foliation into
horizontal (or vertical) lines.

Apart from this technical purpose, circloids are also of an independent
interest, since they may appear as invariant or periodic sets of a toral
homeomorphism. This provides a natural generalisation of the concept of an
invariant essential simple closed curve. Altogether, this leads to the following
Poincar\'e-like classification of conservative pseudo-rotations with bounded
mean motion.
\enlargethispage*{1000pt}
\begin{thmc}
  Suppose $f\in\homtwo$ is a conservative pseudo-rotation with rotation vector
  $\rho \in \R^2$ and bounded mean motion. Then the following hold. 
  \romanlist
\item $\rho$ is totally irrational if and only if $f$ is semi-conjugate to
  $R_\rho$.
\item $\rho$ is neither totally irrational nor rational if and only if $f$ has a
  periodic circloid.
\item $\rho$ is rational if and only if $f$ has a periodic point.  \listend
\end{thmc}
\pagebreak

Finally, the same concepts lead to the following trichotomy for the dynamics of
non-wandering toral homeomorphisms. (We say a map $f$ is {\em non-wandering} if
there exists no non-empty open set $U$ with $f^n(U) \cap U = \emptyset \ \forall
n\geq 1$.)
\begin{thmd}
  \label{t.pp-circloid-transitive} Suppose $f\in\homtwo$ is non-wandering. Then
  one of the following holds.  \romanlist
\item $f$ is topologically transitive;
\item $f$ has two disjoint periodic circloids;
\item $f$ has a periodic point.
\listend
\end{thmd}\noindent
We note that alternatives (i) and (ii) are mutually exclusive, but may both
coexist with (iii). An equivalent way of expressing (ii) is to say that there
exist two disjoint periodic embedded open annuli, both of which contain an
essential simple closed curve.

The existence of a periodic circloid forces the rotation set to be contained in
a line segment which contains no totally irrational rotation vectors (see
Proposition~\ref{p.invariant-circloid} and Remark~\ref{r.invariant-circloid}
below). Hence, we obtain the following corollary.
\begin{core}
  Any non-wandering irrational pseudo-rotation $f\in\homtwo$ is topologically
  transitive.
\end{core}

\noindent {\bf Acknowledgements.} I would like to thank an anonymous referee for
several thoughtful remarks.  This work was supported by a research fellowship
(Ja 1721/1-1) of the German Science Foundation (DFG).

\section{The minimal case}
\label{Minimal} 

The aim of this section is to prove a slightly more general version of
Proposition~A, which also takes into account the situation where the rotation
set is not reduced to a single point, but contained in some lower-dimensional
hyperplane. We define the rotation set of a toral homeomorphism $f\in\homd$,
with lift $F$, on a subset $K\ssq\T^d$ as
\begin{equation}
  \label{eq:3}
   \rho_K(F) \ := \ \left\{ \rho \in \R^d \left|\ \exists n_i \nearrow \infty,\ x_i
  \in K : \iLim (F^{n_i}(x_i) - x_i)/n_i = \rho \right. \right\} \ .
\end{equation}
When $K=\T^d$, this coincides with the standard definition (see
\cite{misiurewicz/ziemian:1989}). Note that for a different lift $F'$ of $f$,
the rotation set $\rho_K(F')$ will be an integer translate of $\rho_K(F)$.
However, this slight ambiguity will not cause any problems, and we will
nevertheless call $\rho_K(F)$ {\em the} rotation set of $f$. Now, suppose
$\rho_K(F)$ is contained in a $d-1$-dimensional hyperplane, that is $\rho_K(F)
\ssq \lambda v + \{v\}^\perp$ for some $v\in\R^d\smin\{0\}$ and $\lambda \in
\R$. In this case, we let
\begin{equation}
  \label{eq:5} D_v(n,z) \ := \ \langle F^n(z)-z-n\rho,v\rangle \ ,
\end{equation}
where $\rho\in\rho_K(F)$ is arbitrary. We say $f$ has {\em bounded mean motion
  parallel to $v$ on $K$} if there exists a constant $c>0$ such that
\begin{equation}
  \label{eq:4}
  |D_v(n,z)| \ \leq \ c \quad \forall n\in\Z,\ z\in K \ .
\end{equation}
By $\|v\|$, we denote the Euclidean norm of a vector $v\in\R^d$, by $\pi_i$ the
projection to the $i$-th coordinate (on any product space). $\pi:\R^d \to \T^d =
\R^d / \Z^d$ will denote the quotient map.

Recall that when $\varphi$ and $\psi$ are endomorphisms of topological spaces
$X$ and $Y$, respectively, then a continuous and onto map $h:X\to Y$ is called a
semi-conjugacy from $\phi$ to $\psi$, if $h\circ \phi = \psi \circ h$. In
general, the existence of a semi-conjugacy from $f_{|K}$ to an irrational
rotation $R_\rho$ does not have any implications for the rotation
set. Therefore, we will use the notion of a {\em regular} semi-conjugacy, which
we define as follows.

Suppose $f \in \homd$ leaves $K\ssq \T^d$ invariant and $R_\rho$ is a rotation
on the $k$-dimensional torus $\T^k$. If $B$ is a $k\times d$ matrix with integer
entries, then a semi-conjugacy $h : K \to \T^k$ from $f_{|K}$ to $R_\rho$ is
called {\em regular with respect to $B$} if it has a lift $H:\pi^{-1}(K) \to
\R^k$ that semi-conjugates $F_{|\pi^{-1}(K)}$ to the translation $T_\rho : z
\mapsto z + \rho$ and satisfies $\sup_{z\in \pi^{-1}(K)} \|H(z) - B(z)\| \leq
\infty$.  Note that in this case $\rho_K(F) \ssq B^{-1}(\rho)$ and $f$ has
bounded mean motion orthogonal to $B^{-1}(\rho)$ (that is, parallel to all $v\in
B^{-1}(\rho)^\perp$).  Furthermore, if $\rho$ is totally irrational, then $B$ is
surjective and hence $B^{-1}(\rho)$ is a $(d-k)$-dimensional hyperplane. When
$B$ is just the projection to the first $k$ coordinates, we simply say that $h$
is regular.

\begin{prop}
  \label{t.minimal-general} Let $f \in \homd$ and $K\ssq\T^d$ be a minimal set
  of $f$. Suppose that there exists an integer vector $v\in\Z^d \smin \{0\}$
  with $\gcd(v_1,\ldots,v_d) =1$ and a number $\rho_0 \in \R\smin\Q$, such that
  \[
     \rho_K(F) \ \ssq \ \frac{\rho_0}{\|v\|^2} \cdot v + \{v\}^\perp \ .
  \]
  Further, assume that $f$ has bounded mean motion parallel to $v$ on $K$. Then
  $f_{|K}$ is regularly semi-conjugate to the one-dimensional irrational
  rotation $r_{\rho_0}: x \mapsto x+\rho_0 \mod 1$.
\end{prop}
The statement can be obtained as a consequence of the Gottschalk-Hedlund Theorem, but
we prefer to give a short direct proof.

\proof First, assume that $v= e^1 = (1,0\ld 0)$. Define $H : K \to \R$ by
\begin{equation}
  \label{eq:6}
  H(z) \ = \ \pi_1(z) + \sup_{n\in\Z} D_{e^1}(n,z) \ = \ \sup_{n\in\Z} 
  \left(\pi_1\circ F^n(z) - n\rho_0\right) \ . 
\end{equation}
Due to the bounded mean motion property $H$ is well-defined, and it is easy to
check that $H\circ F(z) = H(z)+\rho_0$. Furthermore $|H(z) - \pi_1(z)| \leq c$,
where $c$ is the bounded mean motion constant. It remains to show that $H$ is
continuous. In order to do so, note that the function $\varphi(z) =
\sup_{n\in\Z}D_{e^1}(n,z)$ is lower semi-continuous, and $\psi(z) =
\inf_{n\in\Z} D_{e^1}(n,z)$ is upper semi-continuous. Therefore $\varphi-\psi$
is lower semi-continuous, and a straightforward computation shows that it is
furthermore invariant. Since $f_{|K}$ is minimal, this implies that
$\varphi-\psi$ is equal to a constant on $K$, say $c$. It follows that $\varphi
= c+\psi$ is also upper semi-continuous, hence continuous, and thus the same
holds for $H(z) = \pi_1(z) + \varphi(z)$. Since $H$ also satisfies $H(z+v) =
H(z)+\pi_1(v) \ \forall v \in \Z^d$, its projection $h$ to $\T^d$ yields the
required regular semi-conjugacy. The surjectivity of $h$ follows from the
minimality of $r_{\rho_0}$.

In order to reduce the general case to the one treated above, let
$\conv^*(z_1\ld z_n) := \conv(z_1\ld z_n) \smin \{z_1 \ld z_n\}$, where $\conv$
denotes the convex hull. Choose a basis $w^2,\ldots,w^d \in \Z^d$ of
$\{v\}^\perp$ with the property that the $\conv^*(w^2,\ldots,w^d)$
contains no integer vectors.
Next, choose some vector $w^1$, such that $\conv^*(w^1\ld w^d)$ contains no
integer vectors either.
If we denote the matrix $(w^1\ld w^d)$ by $A$, then the latter implies that the
linear toral automorphism $f_A$ induced by $A$ is bijective, such that $\det
A=1$. Furthermore, $\tilde F = A^{-1} \circ F \circ A$ is the lift of a toral
homeomorphism $\tilde f$. There holds
\[
\rho_{f_A^{-1}(K)}(\tilde F) \ = \ A^{-1}(\rho_K(F)) \  \ssq \
\frac{\rho_0}{\|v\|^2} \cdot A^{-1}(v) + \{e^1\}^\perp 
\]
and using $\left(A^{-1}\right)^te^1 \in
\left(A\left(\{e^1\}^\perp\right)\right)^\perp = \R v$ it is easy to check that
$\tilde f$ has bounded deviations parallel to $e^1$.
Thus, it only remains to show that $\langle A^{-1}(v),e^1\rangle = \|v\|^2$. Let
$\tilde v$ be the vector representing the linear functional $x \mapsto
\det(x,w^2\ld w^d)$ on $\R^d$, that is $\det(x,w^2\ld w^d) = \langle x,\tilde v
\rangle \ \forall x \in \R^d$. Then $\tilde v \perp w^i \ \forall i=2\ld d$, and
hence $\tilde v \in \R v$. Furthermore, $x \mapsto \det(x,w^2\ld w^d)$ maps
integer vectors to integers, which implies $\tilde v \in \Z^d$. Finally, the
existence of a vector $w^1 \in \Z^d$ with $\langle w^1,\tilde v\rangle = \det A
= 1$ implies that the coordinates of $\tilde v$ are relatively prime, and hence
$\tilde v = \pm v$. It follows that $|\det(v,w^2\ld w^d)| = \langle v,v\rangle =
\|v\|^2$, and since $\det A=1$ we obtain
\[
    |\langle A^{-1}(v),e^1\rangle| \ = \ |\det(A^{-1}(v),e^2\ld e^d)| \ = \
    |\det(v,w^2\ld w^d)| \ = \ \|v\|^2 \ .
\]
If the sign of $\langle A^{-1} v,e^1\rangle$ is negative, then we simply replace
$w_1$ by $-w_1$. 

Now, as we showed above, there exists a regular semi-conjugacy $h$ from $\tilde
f$ to $r_{\rho_0}$. Thus $h\circ A^{-1}$ yields the required semi-conjugacy from
$f$ to $r_{\rho_0}$, which is regular with respect to $B=\pi_1 \circ A^{-1}$.

\qed

\begin{bem}
  Even without the minimality assumption, the proof of Proposition~A still yields
  the existence of a {\em `measurable semi-conjugacy'}, that is, a measurable
  map $h:K\to\kreis$ that satisfies $h\circ f_{|K} = r_{\rho_0} \circ h$. Since
  $h$ must map any $f_{|K}$-invariant measure $\mu$ to the Lebesgue measure on
  $\T^1$, this is already sufficient to exclude certain exotic behaviour, like
  weak mixing (see \cite{fayad:2002} for examples of this type).
\end{bem}

We obtain the following corollary, which in particular implies Proposition~A.
\begin{cor}
  \label{c.minimal-general} Let $F$ be a lift of $f\in\homd$ and suppose there
  exist vectors $v^1 \ld v^k$ with $\gcd(v^i_1 \ld v^i_d) =1 \ \forall i = 1\ld
  k$ and a totally irrational vector $\rho \in R^{k}$, such that
  \[
  \rho_K(F) \ \ssq \ \bigcap_{i=1}^k \left(\frac{\rho_i}{\|v^i\|^2} \cdot v^i +
  \{v^i\}^\perp\right) \ .
  \]
  Then $f$ is regularly semi-conjugate to the $k$-dimensional irrational
  rotation $R_\rho$.
\end{cor}
\proof Let $h_i$ be the semi-conjugacy between $f$ and $r_{\rho_i}$, obtained
from Proposition~\ref{t.minimal-general} with $v = v^i$. Then $h: \T^d \to \T^k , \
z \mapsto (h_1(z) \ld h_k(z))$ yields the required semi-conjugacy between $f$
and $R_\rho$. Again, the surjectivity of $h$ follows from the minimality of
$R_\rho$, and the regularity is inherited from that of $h_1 \ld h_k$.

\qed
\medskip

The following result is contained in \cite{misiurewicz/ziemian:1991}. 
\begin{thm}[Theorem A in \cite{misiurewicz/ziemian:1991}] \label{t.mz} Let $F$
  be a lift of $f\in \homeo_0(\torus)$ and suppose that $\rho(F)$ has non-empty
  interior. Then given any $\rho \in \inte(\rho(F))$ there exists a minimal set
  $M_\rho$ such that $\rho_{M_\rho}(F) = \{\rho\}$ and $f$ has bounded mean
  motion on $M_\rho$.
\end{thm}
The bounded mean motion property is not explicity stated there, but contained in
the proof (see formula (9)). Together with the preceeding statement, this yields
Corollary B.

\section{Invariant circloids}

In the following, we collect a number of statements about circloids, both on the
open annulus $\A = \kreis \times \R$ and on \torus. These results will be
crucial for the proof of Theorem C in the next and of Theorem D at the end of
this section. Before we start, we want to mention a well-known example, namely
the so-called `pseudo-circle' introduced by Bing~\cite{bing:1948}, which shows
that the structure of a circloid may be much more complicated than that of a
simple closed curve. Later Handel \cite{handel:1982} and
Herman~\cite{herman:1986} showed that the pseudo-circle may appear as an
invariant set of smooth surface diffeomorphisms. Nevertheless, we will see below
that circloids have many `nice' properties, which make them an interesting tool
in the study of toral and annular homeomorphisms.

The definition of a circloid on the annulus is more or less the same as on the
torus. However, for convenience we reformulate it, and introduce some more
terminology.  We say a subset $E \ssq \A$ is an {\em annular continuum}, if it
is compact and connected, and $\A \smin E$ consists of exactly two connected
components which are both unbounded. Note that each of the connected components
will be unbounded in one direction (above or below), and bounded in the other.
We say a subset $C \ssq \A$ is a {\em circloid}, if it is an annular continuum
and does not contain any strictly smaller annular continuum as a subset.

We call a set $E\ssq\A$ {\em essential}, if its complement does not contain any
connected component which is unbounded in both directions. (For compact sets,
this coincides with the usual definition that $E$ is not contained in any
embedded topological disk). Now, suppose that $U \ssq \A$ is bounded from below
and its closure is essential. We will call such a set an {\em upper generating
  set} and define its {\em associated lower component} ${\cal L}(U)$ as the
connected component of $\A \smin \overline{U}$ which is unbounded from below.
Similarly, we call a set $L\ssq \A$ which is bounded from above and has
essential closure a {\em lower generating set}, and define its {\em associated
  upper component} ${\cal U}(L)$ as the connected component of $\A \smin
\overline{L}$ which is unbounded from above. We call an open set ${\cal U}$
(respectively ${\cal L}$) an {\em upper} ({\em lower}) {\em hemisphere}, if
${\cal U} \cup \{+\infty\}$ is bounded from below (${\cal L}\cup\{-\infty\}$ is
bounded from above) and homeomorphic to the open unit disk in $\C$.%
\foot{In order to be absolutely correct, we should say {\em `punctured'}
  hemispheres, but we ignore this for the sake of brevity.}
If $U$, respectively $L$, is connected, then ${\cal L }(U)$, respectively ${\cal
  U}(L)$, is a hemisphere in this latter sense. In order to see this, suppose
$\Gamma$ is a Jordan curve in ${\cal L}(U) \cup \{-\infty\}$. Let $D$ be the
Jordan domain in $\bar\A = \A \cup \{-\infty,+\infty\} \simeq \bar \C$ which is
bounded by $\Gamma$ and does not contain $+\infty$. Since $\overline{U}$ is
connected and essential, $D \cap \overline{U} = \emptyset$. Hence $D$ is
contractible to a point in ${\cal L}(U) \cup \{-\infty\}$. This shows that
${\cal U}(L)$ is simply connected, and the assertion follows from Riemann's
Uniformisation Theorem.

The following remark states a number of elementary properties of the above
objects.
\begin{bem}
  \label{r.exhausive-curves}
\alphlist
\item If $U$ is an upper generating set, then there exist disjoint essential
  simple closed curves $\Gamma_n \ssq {\cal L}(U)$, such that $\bigcup_{n\in\N}
  {\cal L}(\Gamma_n) = {\cal L}(U)$. (For example, the curves $\Gamma_n$ may be
  chosen as the images of the circles with radius $1-1/n$ under the
  homeomorphism from the unit disk to ${\cal L}(U)\cup \{-\infty\}$.) The
  analogous statement holds for lower generating sets.
\item Any annular continuum $E$ is the intersection of a countable nested
  sequence of annuli, bounded by essential simple closed curves. (Simply apply
  (a) to $U=L=E$.)
\item Any upper (lower) hemisphere is an upper (lower) generating set. Hence,
  the expressions ${\cal U}{\cal L}(U)$, ${\cal L}{\cal U}(L)$, ${\cal L}{\cal
    U}{\cal L}(U)$ etc.\ make sense.
\item If $U$ and $U'$ are upper generating sets, then $U' \ssq U$ implies ${\cal
    L}(U) \ssq {\cal L}(U')$. Similarly, if $L$ and $L'$ are lower generating
  sets and $L' \ssq L$, then ${\cal U}(L) \ssq {\cal U}(L')$.
\item If $U$ is an upper separating set, then ${\cal L}(U) \ssq {\cal L}{\cal
    U}{\cal L}(U)$. (Note that ${\cal L}(U) \ssq {\cal U}{\cal L}(U)^c$ by
  definition.)  Similarly, if $L$ is a lower separating set, then ${\cal U}(L)
  \ssq {\cal U}{\cal L}{\cal U}(L)$.
\item Suppose $E$ is both an upper and a lower generating set, for example if
  $E$ is an annular continuum. Then ${\cal L}(E) \ssq {\cal L}{\cal U}(E)$ and
  ${\cal U}(E) \ssq {\cal U}{\cal L}(E)$. (Note that ${\cal L}(E) \ssq
  {\cal U}(E)^c$ and ${\cal U}(E) \ssq {\cal L}(E)^c$.)
  Using (d), this further implies ${\cal U}{\cal L}{\cal U}(E) \ssq {\cal
    U}{\cal L}(E)$ and ${\cal L}{\cal U}{\cal L}(E) \ssq {\cal L}{\cal U}(E)$.
  \listend
\end{bem} 
A general way to obtain circloids is the following. 
\begin{lem}
  \label{l.circloid-generation} Suppose $U$ is an upper generating set. Then
  ${\cal C}^-(U) := \A \smin ({\cal U}{\cal L}(U) \cup {\cal L}{\cal U}{\cal
    L}(U))$ is a circloid. Similarly, if $L$ is a lower generating set, then
  ${\cal C}^+(L) := \A \smin ({\cal L}{\cal U}(L) \cup {\cal U}{\cal L}{\cal
    U}(L))$ is a circloid.
\end{lem}
In particular, every annular continuum $E$ contains a circloid. (Note that
Remarks~\ref{r.exhausive-curves}(e) and (f) imply that $E = \A \smin ({\cal
  U}(E) \cup {\cal L}(E))$ contains both ${\cal C}^+(E)$ and ${\cal C}^-(E)$.)

\proof[Proof of Lemma~\ref{l.circloid-generation}] First, note that since the
operations ${\cal L}$ and ${\cal U}$ always produce hemispheres, ${\cal C}^-(U)$
and ${\cal C}^+(L)$ are annular continua.

Suppose $E$ is an annular continuum which is contained in ${\cal C}^-(U)$. Then,
by definition of ${\cal C}^-(U)$, there holds ${\cal U}{\cal L}(U) \ssq {\cal
  U}(E)$ and ${\cal L}{\cal U}{\cal L}(U) \ssq {\cal L}(E)$. Now ${\cal L}{\cal
  U}{\cal L}(U) \ssq {\cal L}(E)$ implies, due to statement (e) in the preceding
remark, ${\cal L}(U) \ssq {\cal L}(E)$. Hence (d) yields ${\cal U}{\cal L}(E)
\ssq {\cal U}{\cal L}(U)$, and therefore ${\cal U}(E) \ssq {\cal U}{\cal L}(U)$
by (f). Thus ${\cal U}(E) = {\cal U}{\cal L}(U)$.

Similarly, ${\cal U}{\cal L}(U) \ssq {\cal U}(E)$ implies ${\cal L}{\cal U}(E)
\ssq {\cal L}{\cal U}{\cal L}(U)$ by (d) and thus ${\cal L}(E) \ssq {\cal
  L}{\cal U}{\cal L}(U)$ by (f). Hence ${\cal L}(E) = {\cal L}{\cal U}{\cal
  L}(U)$.  Together, we obtain
\[
E \ = \ \A \smin ({\cal U}(E) \cup {\cal L}(E)) \ = 
 \A \smin ({\cal U}{\cal L}(U) \cup  {\cal L}{\cal U}{\cal L}(U)) \ = {\cal C}^-(U) \ .
\]
Of course, the same argument applies to ${\cal C}^+(L)$.

\qed
\medskip

This leads to a nice equivalent characterisation of circloids.  We call
an upper hemisphere $U$ or a lower hemisphere $L$ {\em reflexive}, if ${\cal
  U}{\cal L}(U) = U$ or ${\cal L}{\cal U}(L) = L$, respectively. We call $(U,L)$
a {\em reflexive pair} of hemispheres, if ${\cal U}(L) = U$ and ${\cal L}(U) =
L$. 
\begin{cor}
  \label{l.circloid-reflexive} An annular continuum $C$ is a circloid if and
  only if $({\cal U}(C),{\cal L}(C))$ is a reflexive pair of hemispheres.
\end{cor}

\begin{lem}
  \label{l.thin-continua} Suppose $A$ is an annular continuum with empty
  interior. Then
\[
   {\cal C}^-(A) \ = \ {\cal C}^+(A) \ = \ 
   \partial{\cal U}(A) \cap \partial{\cal L}(A) \ ,
\]
and this is the only circloid contained in $A$. 
\end{lem}
\proof Let $C:= \partial{\cal U}(A) \cap \partial{\cal L}(A)$. Since ${\cal
  U}(A)$ and ${\cal L}(A)$ are open and disjoint, we have
\begin{equation} \label{eq:7} C\ = \ \overline{{\cal U}(A)} \cap \overline{{\cal
      L}(A)} \ = \ \left(\overline{{\cal U}(A)}^c \cup \overline{{\cal
        L}(A)}^c\right)^c \ = \ \A \smin ({\cal LU}(A) \cup {\cal UL}(A)) \ .
\end{equation}
(The last inequality follows from the fact that $\inte(A)=\emptyset$.) We first
show that $C$ is an annular continuum. Since the sets ${\cal LU}(A)$ and ${\cal
  UL}(A)$ are hemispheres, it suffices to prove that their union $V = C^c$ is not
connected. Suppose for a contradiction that it is, and fix two points $z_1 \in
{\cal L}(A) \ssq {\cal LU}(A)$ and $z_2 \in {\cal U}(A) \ssq {\cal
  UL}(A)$. Then, since $V$ is open and connected, we can find an arc $\gamma :
[0,1] \to V$ that joins $z_1$ and $z_2$. However, the sets $\{t\in[0,1] \mid
\gamma(t) \notin \overline{{\cal U}(A)} \}$ and $\{t\in[0,1] \mid \gamma(t)
\notin \overline{{\cal L}(A)} \}$ are both open strict subsets of $[0,1]$ and
their union covers the interval, but they are disjoint (since $\overline{{\cal
    U}(A)} \cup \overline{{\cal L}(A)} = \A$). This contradicts the
connectedness of $[0,1]$. We conclude that $V$ cannot be connected, and hence
$C$ is an annular continuum.

Now ${\cal L}(C) = {\cal LU}(A)$ and ${\cal ULU}(A) \ssq {\cal UL}(A)={\cal
  U}(C)$ by (\ref{eq:7}) and Remark~\ref{r.exhausive-curves}(f). Hence $C \ssq
{\cal C}^+(A)$, and Lemma~\ref{l.circloid-generation} therefore yields $C={\cal
  C}^+(A)$. The same argument shows $C={\cal C}^-(A)$. In particular, $C$ is a
circloid.

Finally, suppose $C'$ is another circloid contained in $A$. Then ${\cal L}(A)
\ssq {\cal L}(C')$, and thus ${\cal L}(A) \cap {\cal U}(C') =
\emptyset$. Therefore
\[ {\cal U}(C') \ \ssq \ {\cal UL}(A)\ = \ {\cal U}(C) \ .
\]
In the same way, we obtain ${\cal L}(C') \ssq {\cal L}(C)$, and hence $C' \ssq
C$. Since $C$ is a circloid, we have $C'=C$.

\qed
\medskip

\noindent
Next, we turn to study circloids which are invariant sets of non-wandering
annular homeomorphisms. Let $\homa$ denote the set of homeomorphisms of $\A$
which are homotopic to the identity. Given $f\in\homa$, an open subset $U\ssq
\A$ is called {\em $f$-wandering}, if $f^n(U) \cap U = \emptyset \ \forall n
\geq 1$. We call $f \in\homa$ {\em non-wandering}, if it does not admit any
non-empty wandering open set, and let $\homar := \{ f\in\homa \mid f \textrm{ is
  non-wandering}\}$. Similarly, we let $\homeo_0^\mathrm{nw}(\torus) := \{
f\in\homtwo \mid f \textrm{ is non-wandering}\}$.  Finally, we call $f\in\homa$
an {\em irrational pseudo-rotation}, if there exists an irrational number
$\rho$, such that for all $z\in\A$ there holds
\begin{equation}
  \label{eq:2a}
  \nLim \pi_1\left(F^n(z) - z\right) /n \ = \ \rho \ .
\end{equation}
Let $p : \R^2 \to \A$ be the canonical projection. The following lemma
will turn out to be useful several times.
\begin{lem}
  \label{l.essentialline} Suppose $f\in\homar$ or $f\in\homeo_0^\mathrm{nw}(\torus)$ has
  no periodic points. Then any open $f$-invariant set contains an essential
  simple closed curve.
\end{lem}
\proof We give the proof for the case of the annulus, the modifications needed
on the torus are minor. Suppose that $f\in\homar$ has no periodic points and $V
\ssq \A$ is an open $f$-invariant set. Fix a small open ball $B \ssq V$. Since
$B$ is non-wandering, there exists some $k\geq 1$ with $f^k(B) \cap B \neq
\emptyset$. Choose a lift $G : \R^2 \to \R^2$ of $f^k$ and a connected component
$\hat B$ of $p^{-1}(B)$, such that $G(\hat B) \cap \hat B \neq \emptyset$. Since
$G$ has no periodic points, a sufficiently small ball $D \ssq \hat B$ will
satisfy $G(D) \cap D = \emptyset$. It follows from a result by Franks
\cite[Prop.~1.3]{franks:1988}, that $G^n(D) \cap D = \emptyset \ \forall
n\in\Z$. Thus, as $p(D)$ is non-wandering for $f^k$, the $G$-orbit of $D$ has to
intersect one of its integer translates. (Note that for any $k\geq 1$, $f$ is
non-wandering if and only if $f^k$ is non-wandering.) The same then certainly
holds for $\hat B$. Since $\bigcup_{n\in\Z} G^n(\hat B) \ssq p^{-1}(V)$ is
connected, this shows that $V$ contains an essential closed curve, which can be
chosen simple.

\qed
\medskip

Since essential simple closed curves are circloids themselves, we obtain the
following corollary.
\begin{cor}
  Suppose $f\in\homar$ or $f\in\homeo_0^\mathrm{nw}(\torus)$ has no periodic points and
  $C$ is an invariant circloid. Then $C$ has empty interior.
\end{cor}
Now we can prove an important property of invariant circloids.
\begin{prop}
  \label{p.circloid-disjointness} Suppose $f\in\homar$ has no periodic points
  and $C_1$ and $C_2$ are $f$-invariant circloids. Then either $C_1=C_2$, or
  $C_1 \cap C_2 = \emptyset$.
\end{prop}
Again, a similar statement holds on the torus, but we will not make use thereof.
\proof First, suppose that ${\cal U}(C_1) \cap {\cal L}(C_2) = {\cal L}(C_1)
\cap {\cal U}(C_2) = \emptyset$. Then ${\cal U}(C_1) \ssq \overline{{\cal
    L}(C_2)}^c$ and therefore ${\cal U}(C_1) \ssq {\cal UL}(C_2) = {\cal
  U}(C_2)$ (the equality comes from Corollary~\ref{l.circloid-reflexive}). In
the same way, we see that ${\cal U}(C_2) \ssq {\cal U}(C_1)$ and thus ${\cal
  U}(C_1) = {\cal U}(C_2)$. The same argument yields ${\cal L}(C_1) = {\cal
  L}(C_2)$, such that $C_1=C_2$.

Otherwise, one of the two intersections is nonempty, we may assume without loss
of generality that $A = {\cal U}(C_1) \cap {\cal L}(C_2) \neq \emptyset$. Since
$A$ is open and invariant, Lemma~\ref{l.essentialline} implies that it contains
an essential simple closed curve $\Gamma$. It is now easy to see that $\Gamma$
separates $C_1$ and $C_2$, that is $C_1 \ssq {\cal L}(\Gamma) $ and $C_2 \ssq
{\cal U}(\Gamma)$, which implies the disjointness of the two sets.

\qed\medskip

In order to apply these results to toral maps, we need the following basic
lemma, whose simple proof is omitted.
\begin{lem}
  \label{l.non-wandering-lift}
  Let $f\in\homeo_0^\mathrm{nw}(\torus)$ and suppose $\rho(F) \ssq \R \times
  \{0\}$ and $f$ has bounded mean motion parallel to $e^2=(0,1)$. Let $\tilde F
  : \A \to \A$ be the (uniquely defined) lift of $f$, such that
  $\sup_{n\in\Z,z\in\A} |\pi_2\circ\tilde F(z)| < \infty$. Then $\tilde F \in
  \homar$.
\end{lem}

We call $f\in\homtwo$ {\em rationally bounded}, if there exists an integer
vector $v$ and some $\lambda \in \Q$, such that $\rho(F) \ssq \lambda v +
\{v\}^\perp$ and $f$ has bounded mean motion parallel to $v$.
\begin{prop} 
  \label{p.invariant-circloid} Suppose $f \in \homeo_0^\mathrm{nw}(\torus)$ has
  no periodic points. Then $f$ is rationally bounded if and only if it has a
  periodic circloid.
\end{prop}
\proof Suppose $f$ is rationally bounded. Using a linear change of coordinates
(as in the proof of Proposition~\ref{t.minimal-general}), we may assume without
loss of generality that $v=e^2$. Suppose $\lambda = p/q$ with $p,q\in\Z$. Let
$\tilde G : \A \to \A$ be the non-wandering lift of $f^q$ provided by
Lemma~\ref{l.non-wandering-lift}~. Then $A:=\bigcup_{n\in\Z} \tilde G^n(\kreis
\times \{0\})$ is invariant, bounded and essential, and thus $C={\cal C}^+(A)$
is an $\tilde F$-invariant circloid. Furthermore,
Proposition~\ref{p.circloid-disjointness} yields $C \cap (C + (0,1)) =
\emptyset$. This implies that there is a simple closed curve $\Gamma$ contained
in the region between $C$ and $C+(0,1)$, whose projection $p(\Gamma)$ will
consequently be contained in $p(C)^c$.  Thus $p(C)$ is the required
$f^q$-invariant circloid.

Conversely, suppose that there exists a $q$-periodic circloid $C$. Then
$\pi^{-1}(C) \ssq \R$ consists of a countable number of connected components,
separated by the lifts of the essential simple closed curve $\Gamma$ contained
in the complement of $C$. A suitable lift $G$ of $f^q$ will leave these
connected components invariant, and it is easy to see that this implies $\rho(G)
\ssq \R v $, where $v \in \Z^2\smin\{0\}$ is the homotopy vector of $\Gamma$.

\qed
\medskip

\begin{bem} \label{r.invariant-circloid} Note that in the above proof, the
  non-existence of periodic points and wandering open sets is only used to
  ensure that the invariant circloid in $\A$ projects down to a circloid in
  $\torus$, via Proposition~\ref{p.circloid-disjointness}. However, this can
  equally be ensured by projecting down only to a sufficiently large finite
  cover of \torus. Hence, even if these assumptions are omitted, we obtain that
  $f\in\homtwo$ is rationally bounded if and only if there exists a lift $\tilde
  f$ of $f$ to a finite cover of $\torus$, such that $\tilde f$ has a periodic
  circloid.
\end{bem}

Theorem~D now follows quite easily from the above results.

\proof[\bf Proof of Theorem D] Suppose that $f\in\homtwo$ has no wandering open
sets. Further, assume that $f$ has no periodic points and is not topologically
transitive. Then there exist two open sets $U_1,U_2$ with disjoint orbit, that
is $\tilde U_1 \cap \tilde U_2 = \emptyset$, where $\tilde U_i =
\bigcup_{n\in\Z} f^n(U_i)$. By Lemma~\ref{l.essentialline}, both $\tilde U_1$
and $\tilde U_2$ contain an essential simple closed curve, which we denote by
$\Gamma_1$ and $\Gamma_2$, respectively. By means of a linear change of
coordinates, we may assume that the homotopy type of these curves is $(1,0)$
(note that since $\Gamma_1$ and $\Gamma_2$ are disjoint, they have the same
homotopy vector). Hence, they lift to essential simple closed curves in
$\A$. Furthermore, any connected component $\hat U_1$ of $\pi^{-1}(\tilde U_1)$
will be contained between two successive lifts of $\Gamma_2$, and consequently
be bounded. A suitable lift $G$ of a suitable iterate of $f$ will leave $\hat
U_1$ invariant. Hence, using Lemma~\ref{l.circloid-generation} we obtain the
existence of two $G$-invariant circloids ${\cal C}^+(\hat U_1)$ and ${\cal
  C}^-(\hat U_1)$. These project to invariant or periodic circloids of $f$. They
cannot project down to the same circloid, because they are both contained in the
region between two successive lifts of $\Gamma_2$.

\qed

\section{The conservative case: Proof of Theorem~C}

Suppose that $f\in\homeo_0(\T^2)$ is a conservative pseudo-rotation with bounded
mean motion. Then the existence of a periodic orbit forces the unique rotation
vector to be rational. Conversely, if the unique rotation vector is rational
then the existence of a periodic orbit follows from a result of
Franks~\cite[Theorem 3.5]{franks:1988a}.  This yields the equivalence in
(iii). (In fact, this holds for pseudo-rotations in general, even without the
conservativity and bounded mean motion hypotheses.) The equivalence in (ii)
follows from Proposition~\ref{p.invariant-circloid} above. Further, if $f$ is
semi-conjugate to a totally irrational rotation on $\T^2$, then the rotation
vector evidently has to be totally irrational. Hence, it remains to prove the
existence of a semi-conjugacy in (i).

Let $\tau : \A \to \torus$ denote the canonical projection and let $T : \A \to
\A,\ (x,y) \mapsto (x,y+1)$. When $A$ is an annular continuum and $B$ is an
arbitrary subset of $\A$, we will use the notation
\begin{eqnarray*}
  A \ \preccurlyeq \ B & :\equi & B\cap  {\cal L}(A) \ = \ \emptyset \ ; \\
  A \ \prec \ B & :\equi & B \ \ssq \ {\cal U}(A) \ .
\end{eqnarray*}
The reverse inequalities are defined analogously. If both $A$ and $B$ are
annular continua and $A \preccurlyeq B$, then we let 
\begin{eqnarray*}
  (A,B) & := & {\cal U}(A) \cap {\cal L}(B) \ ; \\
  {[}A,B{]} & := & \A \smin ({\cal L}(A) \cup {\cal U}(B)) \ .
 \end{eqnarray*}

 (Thus $(A,B)$ is the open region strictly between $A$ and $B$ and $[A,B] =
 (A,B) \cup A \cup B$. One may think of these sets as open and closed
 `intervals' with `endpoints' $A$ and $B$.) Now, suppose $f\in\homtwo$ is an
 irrational pseudo-rotation with rotation vector $\rho$ and bounded mean motion
 with constant $c$. Let $\hat f$ be the lift of $f$ to $\A$ with average
 vertical displacement $\rho_2$, such that $|\pi_2\circ\hat f^n(z)-\pi_2(z)
 -n\rho_2| \leq c \ \forall n\in\Z,z\in\A$. We define
\begin{equation}
A_r \ := \ \bigcup_{n\in\Z} \hat f^n\left(\kreis \times \{r-n\rho_2\}\right) 
\end{equation}
and 
\begin{equation}
C_r \ := \ {\cal C}^+(A_r) \ .
\end{equation}
Note that due to the bounded mean motion property, 
\begin{equation}
  \label{e.circloid-bounds}  A_r \  \ssq \ \kreis \times
  [r-c,r+c] \ .
\end{equation}
Since $A_r$ is also essential, it is a lower generating set, and hence the
definition of $C_r$ makes sense. Further, Lemma~\ref{l.circloid-generation}
implies that the sets $C_r$ are all circloids. The following properties hold and
are easy to verify.
\begin{eqnarray}
  C_{r+1} & = & T(C_r) \label{e.foliation-translates}\\
  \hat f(C_r) & = & C_{r+\rho_2} \label{e.foliation-invariance} \\
  C_r & \preccurlyeq & C_s \quad \textrm{if } r < s \label{e.circloid-disjointness}
\end{eqnarray}
We claim that the circloids $C_r$ are also disjoint, such that 
\begin{equation}
  \label{e.circloid-strict-ieq}
  C_r \ \prec \ C_s \quad \textrm{if } r < s \ .
\end{equation}
This is in fact the crucial point in the proof, and also the part which strongly
relies on the existence of the $f$-invariant measure $\mu$ of full topological
support. In fact, the argument can be seen as a metric version of the one used
in the proof of Lemma~\ref{l.essentialline}. Once we have established this
assertion, the required semi-conjugacy can be constructed quite easily.
\medskip

\noindent {\em Disjointness of the circloids $C_r$.} Note that by going over to
a finite cover of $\torus$ and rescaling, we may assume $c<1/4$. This implies
that $C_r \prec C_{r+1} \ \forall r\in\R$, such that the $C_r$ project down to
circloids on $\torus$. Let $r<s$, and suppose first that $A=[C_r,C_s]$ has empty
interior. In this case Lemma~\ref{l.thin-continua} shows that $A$ contains only
one circloid, and thus $C_r = C_s$. It follows that $C_{r'} = C_{s'} \ \forall
r',s' \in [r,s]$. Choosing $r',s' \in [r,s]$ with $s' = r' + n\rho_2 \bmod 1$ we
obtain $F^n(C_{r'-k}) = C_{r'}$ for some $k\in\Z$. This implies that $f$ has an
invariant or periodic circloid, and is therefore rationally bounded by
Proposition~\ref{p.invariant-circloid}, contradicting the irrationality of
$\rho$.

Thus, we may assume that $A$ has non-empty interior. We claim that $\inte(A)$
contains an essential simple closed curve, which certainly implies the
disjointness of $C_r$ and $C_s$. In order to prove our claim, let $t=(r+s)/2$
and note that, without loss of generality, we may assume $\inte(A') \neq
\emptyset$, where $A'=[C_r,C_t]$ (otherwise, we work with $[C_t,C_s]$; one of
the two sets always has non-empty interior by Baire's Theorem.) Fix some open
ball $V \ssq \inte(A')$ of diameter $\diam(V) \leq 1/8$ and let $V_0 =
\tau(V)$. Choose some integer
\[
  M_1 \ \geq \ \max\left\{\frac{2\mu(\inte(\tau(A)))}{\mu(V_0)},16(c+1)\right\} \ .
\]
Further, choose some integer $m$, such that $(\rho_1',\rho_2') = m\rho \bmod 1$
satisfies $\rho_2' \in \left(0,\frac{t-r}{2M_1^3}\right)$ and $\rho_1' \in
(S\rho_2',2S\rho_2')$, where $$S \ = \ \frac{4M_1(c+1)}{t-r} \ .$$ The fact that
such an $m$ exists follows simply from the minimality of the irrational rotation
$R_\rho$.

Let $G_0:\A \to \A$ be the lift of $f^m$ with rotation vector
$(\rho_1',\rho_2')$, and note that for all $i\leq \frac{t-r}{\rho_2'}$ there
holds $G_0^i(A') = [C_{r+i\rho_2'},C_{t+i\rho_2'}] \ssq A$. Consequently
$f^{im}(\tau(A')) \ssq \tau(A)$, and thus $f^{im}(V_0) \ssq
\inte(\tau(A))$. Since $f^m$ preserves $\mu$, it follows that there exists some
$k\leq M_1$, such that $f^{km}(V_0) \cap V_0 \neq \emptyset$. If
$(\rho_1'',\rho_2'') = k(\rho_1',\rho_2')$, then $\rho_2'' \in
\left(0,\frac{t-r}{2M_1^2}\right)$ and $\rho_1''\in(S\rho_2'',2S\rho_2'')$.

Fix a connected component $\hat V_0 \ssq \R^2$ of $\pi^{-1}(V_0)$ and let
$G_1:\R^2 \to \R^2$ be the lift of $f^{km}$ with $G_1(\hat V_0) \cap \hat V_0
\neq \emptyset$. Then the bounded mean motion property with constant $c\leq
1/4$, which $G_1$ inherits from $f$, implies that $\rho(G_1) =
(\rho_1''',\rho_2''') \in \R^2$ satisfies $\rho_2''' = \rho_2''\in
\left(0,\frac{t-r}{2M_1^2}\right)$ and $\rho_1''' \in
(S\rho_2''',2S\rho_2''')$.%
\foot{It is here were we use $\diam(V) \leq 1/8$ and $M_1 \geq
  16(c+1)$. Otherwise $(\rho_1''',\rho_2''')$ could also be a different lift of
  $(\rho_1'',\rho_2'')$. (Note that $\rho_1''$ is only defined $\bmod 1$.)}
Now, choose $n\in\N$ with $n\in\left[\frac{t-r}{4kM_1\rho_2'}
  ,\frac{t-r}{2kM_1\rho_2'}\right]$. Then $n\rho_2''' \in
\left[\frac{t-r}{4M_1},\frac{t-r}{2M_1}\right]$ and $|\rho_1'''| \geq
Sn\rho_2''' \geq c+1$. The bounded mean motion of $G_1^n$ therefore implies
$G_1^{jn}(\hat V_0) \cap \hat V_0 = \emptyset \ \forall j\in\N$. However, by the
same argument as before there must be some $l\leq M_1$, such that
$f^{lnkm}(V_0)\cap V_0 \neq \emptyset$.  This implies that $G_1^{ln}(\hat V_0)$
has to intersect some integer translate of $\hat V_0$. Since the set $W :=
\bigcup_{j=1}^{ln} G_1^j(\hat V_0)$ is open and connected and $\pi(W) \ssq
\tau(\inte(A))$, this proves our claim.\medskip

\noindent {\em Construction of the semi-conjugacy.} We now define
\begin{equation}
  \label{e.semi-conjugacy}
  H_2(z) \ := \ \sup\{r\in\R \mid z \succ C_r\} \ .
\end{equation}
Using (\ref{e.foliation-translates}) and (\ref{e.foliation-invariance}), it can
easily be checked that
\begin{eqnarray}
  H_2 \circ T(z) & = & H(z) + 1 \ ; \label{e.H-translation} \\
  H_2\circ \hat f(z) & = & H(z) + \rho_2 \ . \label{e.H-semiconjugacy}
\end{eqnarray}
In order to see that $H_2$ is continuous, suppose $(a,b) \ssq \R$ is an open
interval and $z\in H_2^{-1}(a,b)$. Let $c=H_2(z)$, and choose $s\in(a,c)$ and
$t\in(c,b)$. Then $z$ is contained in the open set $(C_s,C_t)$. ($z \succ C_s$
is obvious, and $z \succcurlyeq C_t$ would imply $z \succ C_{t'}$ for all $t' <
t$ by (\ref{e.circloid-strict-ieq}), hence $H_2(z) \geq t$.) However,
$H_2(C_s,C_t) \ssq [s,t] \ssq (a,b)$, such that $H_2^{-1}(a,b)$ contains the
open neighbourhood $(C_s,C_t)$ of $z$. Since $z$ was arbitrary, this proves that
$H_2^{-1}(a,b)$ is open, and as $a,b$ were arbitrary we obtain the continuity
of $H_2$. 

Due to (\ref{e.H-translation}) and (\ref{e.H-semiconjugacy}), $H_2$ projects to
a semi-conjugacy $h_2$ between $f$ and the irrational rotation $r_{\rho_2} : x
\mapsto x+\rho_2 \bmod 1$. In the same way, we can construct a semi-conjugacy
$h_1$ between $f$ and the irrational rotation $r_{\rho_1}$, and $h=(h_1,h_2)$
then yields the required semi-conjugacy between $f$ and $R_\rho$ on $\torus$.

\solidqed

\end{document}